\documentclass[twoside]{article}
\usepackage[dvips]{graphicx}
\usepackage{fancyhdr}

\usepackage{amsmath,amssymb,latexsym}

\usepackage{color}

\usepackage[normalem]{ulem} 

\newtheorem{defn}{\bf Definition}[section]
\newtheorem{lemma}[defn]{\bf Lemma}
\newtheorem{ex}[defn]{\bf Example}
\newtheorem{thm}[defn]{\bf Theorem}
\newtheorem{prop}[defn]{\bf Proposition}
\newtheorem{cor}[defn]{\bf Corollary}

\newcommand{\h}{{\cal H}}
\newcommand{\mn}{\mathbb N}

\def\H{{\cal H}}
\def\h{{\cal H}}
\def\bp{\noindent{\bf Proof: \ }}
\def\ep{\noindent{$\Box$}}
\def\<{\langle}
\def\>{\rangle}

\newcommand{\seqgr[1]}{(#1_n)_{n=1}^\infty}
\def\newin {\,\kern-0.4em\in\kern-0.15em}
\def\newsubset {\kern-0.2em\subset\kern-0.2em}

\def\v{\vspace{.05in}}

\title{On the dual frame induced by an invertible \\ frame multiplier}

\author{Diana T. Stoeva and Peter Balazs \vspace{.1cm} \\ 
  Acoustics Research Institute, Austrian Academy of Sciences,\\
Wohllebengasse 12-14, Vienna A-1040, Austria \\
dstoeva@kfs.oeaw.ac.at, peter.balazs@oeaw.ac.at \\
}

\begin{document}


\maketitle
\thispagestyle{empty}

\markboth{\footnotesize \rm \hfill D. T. STOEVA AND P. BALAZS \hfill}
{\footnotesize \rm \hfill ON THE DUAL FRAME INDUCED BY AN INVERTIBLE FRAME MULTIPLIER \hfill}

\begin{abstract}
Recently it has been established that given an invertible frame multiplier with semi-normalized symbol, a specific dual of any of the two  involved frames can be determined  for the inversion purpose. 
The inverse can be represented as a multiplier with the reciprocal symbol, this particular dual of one of the given frames, and any dual of the other frame. The specific dual is the only one having this property among all Bessel sequences.

In this manuscript we extend the results showing that the specific dual with the above mentioned property is unique among all possible sequences. Furthermore,  the symbol is allowed to be not necessarily semi-normalized. Finally we characterize cases when the canonical dual frame and the new specific dual frame coincide.
\vspace{5mm} 
\vspace{3mm}
\end{abstract}

\section{Introduction}

For a representation of functions, beside orthonormal bases, other options like  Riesz bases \cite{bari51} or frames   \cite{ds52} exist. 
Riesz bases extend orthonormal bases giving up with the orthonormality restrictions, but still providing unique representations. 
Frames extend Riesz bases giving up with the basis-property and allow redundant representations.
The redundancy is very desirable in applications, for example in signal processing \cite{chai10}, in particular for 
noise reduction \cite{BHnoisereduction} or signal modification \cite{mathlaw03}. 
A natural way to modify functions is first to transform them (via some frame), then multiply the resulting coefficients with a fixed scalar sequence, and then transform back into the original domain via the same or another frame. 
The operators used to do such a modification are called multipliers. They are of big importance both in mathematics \cite{Arias2008581,feinow1,xxlbayasg11} and applications \cite{xxllabmask1,wanbro06,majxxl10,oltokr13}. 
The invertibility of such operators is also of mathematical and  practical interest. 
In this paper we deal with the question on how to represent the inverse as a multiplier involving appropriate dual frames of the given ones. 

Riesz multipliers with semi normalized symbols are always invertible and the inverse can be written as a multiplier with the reciprocal symbol and the canonical dual frames \cite{xxlmult1}. Considering the more general case of frame multipliers with semi-normalized symbols, invertibility is not always the case and even in the case of invertibility, the above mentioned way of inversion does not necessarily work, though it is still possible in certain cases. Investigation of the frame multiplier case was done in \cite{bsreprinv2015} where formulas for the inverse were determined involving appropriate dual frames of the given ones. Uniqueness of these dual frames was established among the set of the Bessel sequences.   
A step further, the SampTA-proceedings \cite{sbsampta15} was devoted to the investigation of frame multipliers where the symbol is not necessarily semi-normalized and 
non-Bessel sequences were also involved in the consideration. 
The results there were announced without proofs. 
This paper is an extension of \cite{sbsampta15}. 
Here we give proofs of the results, as well as provide new results related to characterization of frame multiplier cases where the inversion can be done using precisely the formula from the Riesz multiplier case.

The paper is organized as follows. 
Section \ref{motiv} presents the existing and motivating results regarding the particular unique duals involved in representing the inverse multiplier.
These results are related to semi-normalized symbols and the assumption of Bessel sequences.   
The section also points to the main questions, concerning generalizations of the existing results, which are treated in the later sections.
In Section \ref{snm} the assumption of Bessel sequences is lifted.
In Section \ref{snnnm} not necessarily semi-normalized symbols are considered. 
And finally in Section \ref{struct0} we comment on and characterize cases when the formula for representing the inverse of a Riesz multiplier, using the canonical duals,   hold for redundant frame multipliers.

\subsection*{Basic Definitions, Notation, and Known Facts}

Throughout the paper, $\h$ denotes a separable Hilbert space,  $(e_n)_{n=1}^\infty$ - an orthonormal basis of $\h$, $\Phi$ (resp. $\Psi$) - a sequence $(\phi_n)_{n=1}^\infty$ (resp. $(\psi_n)_{n=1}^\infty$) with elements from $\h$, $m$ - a complex scalar sequence $(m_n)_{n=1}^\infty$, $\overline{m}$ - the sequence  $(\overline{m_n})_{n=1}^\infty$ consisting of the complex conjugates of the elements of $m$, 
$1/m$ - the sequence $(1/m_n)_{n=1}^\infty$, $m\Phi$ - the sequence $(m_n\phi_n)_{n=1}^\infty$. 

The sequence $m$ is called \emph{semi-normalized} if there exist $a,b$ so that $0<a\leq |m_n|\leq b<\infty$ for all $n\in\mn$. 
The term \emph{operator} is used with the meaning of a linear mapping; an operator $L:\h\to\h$ is called \emph{invertible} if it is bounded and bijective.  

For extensive presentation on frames we refer to \cite{Casaz1,ole1,gr01,he98-1}. 
Recall some needed definitions and facts from frame theory.

The sequence $\Phi$ is called 

\begin{itemize}
\item[-] \vspace{-3.3mm} a \emph{Bessel sequence in $\h$} if there exists $B\in (0,\infty)$ so that $\sum_{n=1}^\infty |\<f, \phi_n\>|^2 \leq B \|f\|^2$ for all $f\in\h$; 
\item[-] \vspace{-3mm} a \emph{frame for $\h$} if it is a Bessel sequence in $\h$ and there exists $A\in (0,\infty)$ so that  $ A \|f\|^2 \leq \sum_{n=1}^\infty |\<f, \phi_n\>|^2$ for all $f\in\h$; 
\item[-] \vspace{-3mm}  a \emph{Riesz basis for $\h$} if it is both a frame for $\h$ and a Schauder basis for $\h$. 
\end{itemize}
\vspace{-2mm} 
Two frames $\Phi$ and $\Psi$ for $\h$ are called \emph{equivalent} when there exists an invertible operator $L:\h\to\h$ mapping $\phi_n$ to $\psi_n$, $n\in\mn$. 

A  frame $(f_n)_{n=1}^\infty$ which satisfies $f=\sum_{n=1}^\infty \<f, \phi_n\>f_n =  \sum_{n=1}^\infty \<f, f_n\>\phi_n$ for all $f\in\h$ and a given frame $\Phi$ is called a \emph{dual frame of $\Phi$}. The \emph{frame operator of a frame $\Phi$}, denoted by $S_\Phi$, is determined by $S_\Phi f := \sum_{n=1}^\infty \<f,\phi_n\>\phi_n$; it is well defined on $\h$ and invertible \cite{ds52}. The sequence 
$\widetilde{\Phi}:=(S_\Phi^{-1}(\phi_n))_{n=1}^\infty$ is a dual frame of $\Phi$ \cite{ds52}, called the \emph{canonical dual of $\Phi$}.

In addition to dual frames, we also take into consideration sequences which are not necessarily frames and still provide series expansions (examples can be found e.g. in \cite{LO,Scharact,S2009IJPAM,bstable09}). Given a frame $\Phi$ for $\h$, the sequence $F=(f_n)_{n=1}^\infty$ with elements from $\h$ is called 

- an \emph{analysis pseudo-dual} (in short, \emph{a-pseudo-dual}) \emph{of $\Phi$}, if  $f= \sum_{n=1}^\infty \<f, f_n\> \phi_n$ for all $f\in\h$.

- a \emph{synthesis pseudo-dual} (in short, \emph{s-pseudo-dual}) \emph{of $\Phi$}, if  $f= \sum_{n=1}^\infty \<f, \phi_n\> f_n$ for all $f\in\h$.

For more on a-pseudo-duals and s-pseudo-duals, see \cite{Scharact}.

Given sequences $m$, $\Phi$, and $\Psi$, the operator $M_{m,\Phi,\Psi}$ determined by 
$$M_{m,\Phi,\Psi} f = \sum_{n=1}^\infty m_n\<f,\psi_n\>\phi_n, \ f\in {\rm dom}(M_{m,\Phi,\Psi})\subseteq \h$$
is called a \emph{multiplier} and $m$ is called the \emph{symbol} of $M_{m,\Phi,\Psi}$. The multiplier $M_{m,\Phi,\Psi}$ is called \emph{well-defined on $\h$} if its domain is $\h$.  
When $\Phi$ and $\Psi$ are frames (resp. Riesz bases) for $\h$,  $M_{m,\Phi,\Psi}$  is called a \emph{frame multiplier}, resp. \emph{Riesz multiplier}. Frame multipliers with semi-normalized symbols are well-defined on $\h$ \cite{xxlmult1}. When $m$ is the constant sequence $(c,c,c,\ldots)$, the multiplier $M_{m,\Phi,\Psi}$ will be denoted by $M_{(c),\Phi,\Psi}$.  For works oriented to invertibility of multipliers refer to \cite{bsreprinv2015,balsto09new,uncconv2011,bstable09}.

\sloppy
Notice that the assumptions of all assertions in the next sections lead to well-defindedness of the multipliers $M_{m,\Phi,\Psi}$ and we will not mention this explicitly in the statements. 

\section{Motivation - the case of semi-normalized symbols and Bessel sequences}\label{motiv}

Here we review the motivating results which concern multipliers with semi-normalized symbols and take into consideration only Bessel sequences.  
 The starting point for considering the topic was a result about invertibility of multipliers for Riesz bases and representation of the inverse via the canonical duals:
 
 \begin{prop} \label{rb} {\rm \cite{xxlmult1}}
 Let $\Phi$ and $\Psi$ be Riesz bases for $\h$, and let 
the symbol $m$ be semi-normalized.  
Then $M_{m, \Phi, \Psi}$ is invertible and 
\begin{equation}\label{minvrb}
M_{m, \Phi, \Psi}^{-1}  = M_{1/m, \widetilde{\Psi}, \widetilde{\Phi}}.
\end{equation} 
 \end{prop}
 
 The above statement has lead to investigation of multipliers for (possibly overcomplete) frames aiming representation of the inverse (in case of invertibility) of a similar type like (\ref{minvrb}), and
 has led to the following results:

\begin{thm} \label{ff} {\rm \cite{bsreprinv2015}}
Let $\Phi$ and $\Psi$ be frames for $\h$, and let 
the symbol  
$m$ be semi-normalized.  
Assume that $M_{m, \Phi, \Psi}$ is invertible. 
Then the following statements hold.
\begin{itemize}
\item[{\rm (a1)}]
There exists a unique dual frame $\Phi^\dagger$ of $\Phi$, 
so that 
\begin{equation}\label{minv1}
M_{m, \Phi, \Psi}^{-1}  = M_{1/m, \Psi^d, \Phi^\dagger}, \ \forall \mbox{ dual frames $\Psi^d$ of $\Psi$}.
\end{equation} 
The frame $\Phi^\dagger$ equals
the sequence $((M_{m,\Phi,\Psi}^{-1})^*(\overline{m_n}\psi_n))_{n=1}^\infty$.
\item[{\rm (b1)}]
There exists a unique dual frame $\Psi^\dagger$ of $\Psi$, 
so that  
\begin{equation}\label{minv2}
M_{m, \Phi, \Psi}^{-1}  = M_{1/m, \Psi^\dagger, \Phi^d}, \ \forall \mbox{ dual frames $\Phi^d$  of $\Phi$}.
\end{equation} 
The frame $\Psi^\dagger$ equals the sequence $(M_{m,\Phi,\Psi}^{-1}(m_n\phi_n))_{n=1}^\infty$.  
\end{itemize}

\begin{itemize}
\item[{\rm (a2)}]
 $\Phi^\dagger$ is the only Bessel sequence in $\h$ which satisfies 
(\ref{minv1}).
\item[{\rm (b2)}]
 $\Psi^\dagger$ is the only Bessel sequence in $\h$ which satisfies 
(\ref{minv2}).
\item[{\rm (a3)}]
If $F=(f_n)_{n=1}^\infty$ is a Bessel sequence in $\h$ which fulfills 
 $$M_{m, \Phi, \Psi}^{-1}  = M_{1/m, F, \Phi^\dagger},$$
then $F$ must be a dual frame of  $\Psi$. 
\item[{\rm (b3)}]
If $G=(g_n)_{n=1}^\infty$ is a Bessel sequence in $\h$ which fulfills 
$$M_{m, \Phi, \Psi}^{-1}  = M_{1/m, \Psi^\dagger, G},$$
then $G$ must be a dual frame of $\Phi$.
\end{itemize}\end{thm}

 \vspace{.05in}
 Theorem \ref{ff} naturally leads to the following questions:
 
\begin{quote} 
[{$\bf Q1$}] {\em 
Are there non-Bessel sequences, which can be used in the role of $\Phi^d$ and $\Psi^d$, or in the role of $\Phi^\dagger$ and $\Psi^\dagger$?
 } 
\end{quote}
\begin{quote} 
 [{$\bf Q2$}] {\em 
What about cases when  $m$ is not necessarily semi-normalized? 
} 
\end{quote}

Under the assumptions of Theorem \ref{ff}, the uniquely determined dual frame $\Phi^\dagger$ of $\Phi$ (resp. $\Psi^\dagger$ of $\Psi$) is referred to as \emph{the dual frame of $\Phi$} (resp. \emph{$\Psi$}) \emph{induced by $M_{m,\Phi,\Psi}$}.  

While the representations of the inverse in (\ref{minv1}) and (\ref{minv2}) determine an appropriate analogue of (\ref{minvrb}) for the general frame case, 
notice that one may still have precisely the formula (\ref{minvrb}) even in the overcomplete case (e.g., for every frame $\Phi$ for $\h$,
 $M_{(1),\Phi,\Phi}$ is the frame operator of $\Phi$ and so it is invertible and fulfills (\ref{minvrb})), but this is certainly not always the case, see \cite[Example 4.2]{bsreprinv2015}. 
The statement below presents some sufficient conditions, in certain cases being also necessary, for validity of   (\ref{minvrb}).

\v
\begin{prop} \label{q} {\rm \cite{bsreprinv2015}}
Let the assumptions of Theorem \ref{ff} hold.  
Then 
\begin{equation} \label{psid}
M_{m,\Phi,\Psi}^{-1}=M_{1/m,\widetilde{\Psi},\widetilde{\Phi}} \Leftarrow 
 \Psi \mbox{ is equivalent   to } m\Phi 
\Leftrightarrow
 \Psi^\dagger = \widetilde{\Psi};
\end{equation}
\begin{equation} \label{phid}
M_{m,\Phi,\Psi}^{-1}=M_{1/m,\widetilde{\Psi},\widetilde{\Phi}} \Leftarrow 
\Phi \mbox{ is equivalent   to } \overline{m}\Psi
\Leftrightarrow 
\Phi^\dagger = \widetilde{\Phi}.
\end{equation}

\vspace{.07in}
\noindent 
For the case $m=(c,c,c,\ldots)$, $c\neq 0$,  
\begin{eqnarray*}
M_{(c),\Phi,\Psi}^{-1}=M_{(1/c),\widetilde{\Psi},\widetilde{\Phi}} \Leftrightarrow \Psi \mbox{ is equivalent   to } \Phi
&\Leftrightarrow&  \Psi^\dagger = \widetilde{\Psi} 
\\
&\Leftrightarrow& \Phi^\dagger = \widetilde{\Phi}.
\end{eqnarray*} 
\end{prop}

Notice that the  equivalence of the frames $\Psi$ and $m\Phi$ (resp. $\Phi$ and $\overline{m}\Psi$) already implies invertibility of $M_{m,\Phi,\Psi}$, 
 and
(\ref{minvrb}) holds. 
When $\Psi$ is not equivalent to $m\Phi$ and $\Phi$ is not equivalent to $\overline{m}\Psi$, then $M_{m,\Phi,\Psi}$ can still be invertible (see, e.g.,  
\cite[Ex. 4.6.3(ii)]{bstable09}, where neither $\Psi$ is equivalent to $m\Phi$, nor $\Phi$ is  equivalent   to $\overline{m}\Psi$, but $M_{m,\Phi,\Psi}$ is invertible; furthermore (\ref{minvrb}) does not hold) and the question is whether in such cases it may happen at all to have (\ref{minvrb}):

\begin{quote} 
[{$\bf Q3$}] {\em 
Are the first implications in (\ref{psid}), resp. (\ref{phid}), actually equivalences?
 } 
\end{quote}

\vspace{.1in} As shown in Proposition \ref{q}, in the case when $m$ is a constant sequence, the answer of $[\bf Q3]$ is affirmative, so the aim now is to answer the question in the general case of $m$ not necessarily being a constant sequence.

The paper is devoted to the above and related questions. 
Answers of $[Q1]$ and $[Q2]$ were announced in \cite{sbsampta15} without proofs; they are presented here. In addition, in this paper we also present investigation related to the question $[Q3]$.

\section{On [Q1]: non-Bessel sequences} \label{snm}

In this section we extend Theorem \ref{ff}, taking non-Bessel sequences into consideration, and answer   $[Q1]$, while continuing the numbering from Theorem \ref{ff}.
 
 \vspace{.05in}
 \begin{thm} \label{ff2} 
Let the assumptions of Theorem \ref{ff} hold and let $\Phi^\dagger$ and $\Psi^\dagger$ be determined by Theorem \ref{ff}.
Then the following statements hold.

\begin{itemize}
\item[{\rm (a4)}]
 $\Phi^\dagger$ is the only  sequence in $\h$ which satisfies 
(\ref{minv1}).
\item[{\rm (b4)}]
 $\Psi^\dagger$ is the only sequence in $\h$ which satisfies 
(\ref{minv2}).
\end{itemize}

\begin{itemize}
\item[{\rm (a5)}]
$M_{m, \Phi, \Psi}^{-1}  = M_{1/m, \Psi^{sd}, \Phi^\dagger}, \ \forall \mbox{ s-pseudo-duals $\Psi^{sd}$ of $\Psi$}$.
\item[{\rm (b5)}]
$M_{m, \Phi, \Psi}^{-1}  = M_{1/m, \Psi^\dagger, \Phi^{ad}}, \ \forall \mbox{ a-pseudo-duals $\Phi^{ad}$  of $\Phi$}.$
\end{itemize}

\begin{itemize}
\item[{\rm (a6)}]
If $F=(f_n)_{n=1}^\infty$ is a sequence in $\h$ such that $M_{1/m, F, \Phi^\dagger}$ is well-defined and
 $$M_{m, \Phi, \Psi}^{-1}  = M_{1/m, F, \Phi^\dagger},$$
then $F$ must be an s-pseudo-dual of  $\Psi$.
\item[{\rm (b6)}]
If $G=(g_n)_{n=1}^\infty$ is a sequence in $\h$ such that $ M_{1/m, \Psi^\dagger, G}$ is well-defined and  
$$M_{m, \Phi, \Psi}^{-1}  = M_{1/m, \Psi^\dagger, G},$$
then $G$ must be an a-pseudo-dual of $\Phi$.
\end{itemize}
\end{thm} 

For the proof of the above proposition, we will use the following statements:

\begin{lemma} \label{lemnullop}
Let $\Phi$ be a frame for $\h$ and let $m$ be such that $m_n\neq 0$ for every $n$.
 Assume that $G=(g_n)_{n=1}^\infty$ is a sequence in $\h$ so that $M_{m, G, \Phi^d}$ (resp.  $M_{m,  \Phi^d, G}$) is the null-operator on $\h$ for every dual frame $\Phi^d$ of $\Phi$. Then $G$ is the null-sequence. 
\end{lemma} 
\bp  Assume that the sequence $G$ is such that  $M_{m, G, \Phi^d}$ is the null operator on $\h$ for every dual frame $\Phi^d=(\phi_n^d)_{n=1}^\infty$ of $\Phi$.  
Let $f\in\h$. Then  
\begin{equation}\label{a1}
\sum_{n=1}^\infty m_n \<f,\phi_n^d\>g_n = 0, \ \forall \mbox{ dual frame $\Phi^d$ of $\Phi$}.
\end{equation}
The dual frames of $\Phi$ are characterized as the sequences $(\widetilde{\phi}_n + h_n - \sum_{j=1}^\infty \<\widetilde{\phi}_n, \phi_j\>h_j)_{n=1}^\infty$ where $(h_n)_{n=1}^\infty$ runs through the Bessel sequences (see, e.g., \cite[Theor. 5.6.5]{ole1}). 
Using this characterization, we get 
$$
\sum_{n=1}^\infty m_n \<f,\widetilde{\phi}_n + h_n - \sum_{j=1}^\infty \<\widetilde{\phi}_n, \phi_j\>h_j\> g_n = 0, \ \forall 
\mbox{ Bessel sequence $(h_n)_{n=1}^\infty$ in $\h$}.
$$
Using (\ref{a1}) with $\widetilde{\Phi}$ in the role of $\Phi^d$, it then follows that 
$$
\sum_{n=1}^\infty m_n \<f, h_n - \sum_{j=1}^\infty \<\widetilde{\phi}_n, \phi_j\>h_j\> g_n = 0, \ \forall 
\mbox{ Bessel sequence $(h_n)_{n=1}^\infty$ in $\h$}.
$$
Take $(h_n)_{n=1}^\infty=(e_1, 0, 0, 0, \ldots,)$. 
Then
$$
m_1 \<f, e_1\> g_1 - 
\sum_{n=1}^\infty m_n \<f, \<\widetilde{\phi}_n, \phi_1\>e_1\> g_n = 0, 
$$
implying 
$$
m_1 \<f, e_1\> g_1 - 
 \<f, e_1\> \sum_{n=1}^\infty m_n \<\phi_1, \widetilde{\phi}_n\> g_n = 0, 
$$
which by (\ref{a1}) leads to  
$m_1 \<f, e_1\> g_1=0.$
Choosing $f=e_1$, one comes to the conclusion that $g_1=0$.

In a similar way, taking $(h_n)_{n=1}^\infty= (0, 0, \ldots, 0, e_i, 0, \ldots)$ with $e_i$ being at the $i$-th place, $i\geq 2$, we get
 $g_i=0$, $i\geq 2$.

Assuming that $G$ is a sequence such that  $M_{m, \Phi^d, G}$ is the null operator on $\h$ for every dual frame $\Phi^d=(\phi_n^d)_{n=1}^\infty$ of $\Phi$, one may proceed with similar techniques to prove that $G$ is the null sequence.
\ep

 \vspace{.1in} \noindent {\bf Proof of Theorem \ref{ff2}: }
Throughout the proof, let $M$ denote the multiplier $M_{m,\Phi,\Psi}$.

  \v
(a4) 
Assume that $F=\seqgr[f]$ is a sequence in $\h$ so that 
$M^{-1}  = M_{1/m, \Psi^d, F}$ for every dual frame $\Psi^d$ of $\Psi$. 
By Theorem \ref{ff} (b1), it follows that 
$M_{1/m, \Psi^\dagger, \Phi^d}= M_{1/m, F, \Phi^d}$ for every dual frame $\Psi^d$ of $\Psi$. Then Lemma \ref{lemnullop} implies that $F=\Psi^\dagger$.

\v 
(b4) can be proved in a similar way as (a4).

\sloppy

\v
(a5) Assume that $\Psi^{sd}$ is an s-pseudo-dual of $\Psi$. 
Then for every $f\in\h$, $\sum_{n=1}^\infty \<f, \psi_n\> \psi_n^{sd}=f$, which implies that
$$M^{-1}f  =\sum_{n=1}^\infty \<M^{-1}f, \psi_n\> \psi_n^{sd} 
=  \sum_{n=1}^\infty \frac{1}{m_n} \<f, (M^{-1})^*(\overline{m_n}  \psi_n)\> \psi_n^{sd}.$$
Since $\Phi^\dagger$ is  the sequence $( (M^{-1})^*(\overline{m_n}  \psi_n))_{n=1}^\infty$, it now
follows that $ M_{1/m, \Psi^{sd}, \Phi^\dagger}$ is well defined and equal to $M^{-1}$.

\v
(b5) Assume that $\Phi^{ad}$ is an a-pseudo-dual of $\Phi$. 
Then for every $f\in\h$, $\sum_{n=1}^\infty \<f, \phi^{ad}_n\> \phi_n=f$, which implies that
$$M^{-1}f 
= \sum_{n=1}^\infty \frac{1}{m_n}\<f, \phi^{ad}_n\> M^{-1}(m_n\phi_n).$$
Since $\Psi^\dagger$ is  the sequence $( M^{-1}(m_n\phi_n))_{n=1}^\infty$, it follows  that $M_{1/m,\Psi^\dagger, \Phi^{ad}}$ is well defined and 
equal to $M^{-1}$.

\v (a6) Assume that $F=(f_n)_{n=1}^\infty$ is a sequence in $\h$ such that $M_{1/m, F, \Phi^\dagger}$ is well-defined and
 $M^{-1}  = M_{1/m, F, \Phi^\dagger}.$
  Then for every $f\in\h$, 
$$M^{-1}f 
= \sum_{n=1}^\infty \frac{1}{m_n} \<f, \phi_n^\dagger\> f_n
= \sum_{n=1}^\infty \frac{1}{m_n} \<f, (M^{-1})^*(\overline{m_n}  \psi_n)\> f_n 
=\sum_{n=1}^\infty  \<M^{-1}f,   \psi_n\> f_n.
$$
Therefore, 
 $F$ is an s-pseudo-dual of  $\Psi$. 

\v (b6) Assume that $G=(g_n)_{n=1}^\infty$ is a sequence in $\h$ such that $ M_{1/m, \Psi^\dagger, G}$ is well-defined and  
$M^{-1}  = M_{1/m, \Psi^\dagger, G}.$ 
Then for every $f\in\h$, 
$$f 
= M M_{1/m, \Psi^\dagger, G} f 
= M \sum_{n=1}^\infty \frac{1}{m_n} \<f, g_n\> M^{-1}(m_n\phi_n)
=   \sum_{n=1}^\infty  \<f, g_n\> \phi_n,$$
which implies that 
 $G$ is an a-pseudo-dual of $\Phi$. 
 \ep

\section{On [Q2]: non-semi-normalized symbols} \label{snnnm}
 
Here we consider multipliers with symbols not necessarily being semi-normalized, as well as taking into consideration sequences not necessarily being Bessel.

\begin{prop} \label{ff3} 
Let $\Phi$ and $\Psi$ be frames for $\h$, and let 
the symbol  
$m$ be such that $m_n\neq 0$ for every $n$ and the sequence $m\Phi$ be a frame for $\h$. 
Assume that $M_{m, \Phi, \Psi}$ is invertible. 
Then the following statements hold.
\begin{itemize}
\item[{\rm (i)}]
There exists a unique sequence $\Psi^\dagger$ in $\h$ 
so that 
\begin{equation*}
 M_{m, \Phi, \Psi}^{-1}  = M_{1/m, \Psi^\dagger, \Phi^{ad}}, \ \forall \mbox{ a-pseudo-duals $\Phi^{ad}$  of $\Phi$,}
\end{equation*}
 and it is a dual frame of $\Psi$. Furthermore, $\Psi^\dagger=(M_{m,\Phi,\Psi}^{-1}(m_n\phi_n))_{n=1}^\infty$.
 
 \item[{\rm (ii)}]
 If $G=(g_n)_{n=1}^\infty$ is a sequence in $\h$  such that $ M_{1/m, \Psi^\dagger, G}$ is well-defined and  
$$M_{m, \Phi, \Psi}^{-1}  = M_{1/m, \Psi^\dagger, G},$$
then $G$ must be an a-pseudo-dual  of $\Phi$.
\end{itemize}
\end{prop} 
\bp
(i) Consider $M_{m,\Phi,\Psi}=M_{(1), m\Phi,\Psi}$. Then Theorem \ref{ff} and Theorem \ref{ff2} can be applied to $M_{(1), m\Phi,\Psi}$, implying that 
there exists a dual frame $\Psi^\dagger$ of $\Psi$ 
so that
\begin{equation}\label{ad1}
M_{(1), m\Phi,\Psi}^{-1} = M_{(1), \Psi^\dagger,(m\Phi)^{ad}}, \ \forall \mbox{ a-pseudo-dual $(m\Phi)^{ad}$ of $m\Phi$}
\end{equation} \sloppy
and this frame $\Psi^\dagger$ is the only sequence in $\h$ satisfying (\ref{ad1}); furthermore, $\Psi^\dagger= (M_{(1), m\Phi,\Psi}^{-1} (m_n\phi_n))_{n=1}^\infty=( M_{m, \Phi,\Psi}^{-1}(m_n\phi_n) )_{n=1}^\infty$.
Since the a-pseudo-duals of $m\Phi$ are precisely the sequences $(1/\overline{m})\Phi^{ad}$ with $\Phi^{ad}$ being an a-pseudo-dual of $\Phi$, the property (\ref{ad1}) is equivalent to
\begin{equation*}
M_{m, \Phi,\Psi}^{-1} = M_{1/m, \Psi^\dagger,\Phi^{ad}}, \ \forall \mbox{ a-pseudo-duals $\Phi^{ad}$ of $\Phi$}.
\end{equation*}

(ii)
Assume that $G=(g_n)_{n=1}^\infty$ is a sequence in $\h$  such that $ M_{1/m, \Psi^\dagger, G}$ is well-defined and  
$M_{m, \Phi, \Psi}^{-1}  = M_{1/m, \Psi^\dagger, G}$,  i.e.,  
$M_{(1), m\Phi, \Psi}^{-1}  = M_{(1), \Psi^\dagger, (1/\overline{m})G}.$ 
By Theorem \ref{ff2}  applied to $M_{(1), m\Phi,\Psi}$, 
we can conclude that $(1/\overline{m})G$ is an a-pseudo-duals of $m\Phi$, 
which implies that $G$ is an a-pseudo-dual  of $\Phi$. 
\ep

\v \noindent {\bf Remark:} In the above proposition it is assumed that $\Phi$ is a frame for $\H$ and $m$ is such that $m\Phi$ is also a frame for $\h$. Notice that by \cite{stoevaxxl09}, both $\Phi$ and $m\Phi$ being frames may happen only in  cases determined as follows:

$\bullet$  $\Phi$ is a Riesz basis for $\h$ and $m$ is semi-normalized (in this case $m\Phi$ is a Riesz basis for $\h$);

$\bullet$  $\Phi$ is an overcomplete frame for $\h$ and $m\in\ell^\infty$, or,  
 $\Phi$ is an overcomplete frame for $\h$ so that no positive constant $a$ fulfills $a<\|\phi_n\|$ for all $n\in\mn$ and $m\notin\ell^\infty$ (in such cases $m\Phi$ can be an overcomplete frame for $\h$).

 \v 
The above Proposition \ref{ff3} is related to the unique dual frame of $\Psi$ induced by $M_{m,\Phi,\Psi}$.
With a similar proof, one can show the validity of the following statement which concerns the unique dual frame of $\Phi$ induced by $M_{m,\Phi,\Psi}$:

\begin{prop} \label{ff32} 
Let $\Phi$ and $\Psi$ be frames for $\h$, and let 
the symbol  
$m$ be such that $m_n\neq 0$ for every $n$ and the sequence $m\Psi$ is a frame for $\h$. 
Assume that $M_{m, \Phi, \Psi}$ is invertible. 
Then the following statements hold.
\begin{itemize}
\item[{\rm (i)}]
There exists a unique sequence $\Phi^\dagger$ in $\h$ 
so that 
\begin{equation*}
 M_{m, \Phi, \Psi}^{-1}  = M_{1/m, \Psi^{sd}, \Phi^\dagger}, \ \forall \mbox{ s-pseudo-duals $\Psi^{sd}$ of $\Psi$,}
\end{equation*}
 and it is a dual frame of $\Phi$.  Furthermore, $\Phi^\dagger=((M^{-1})^*(\overline{m_n}\psi_n))_{n=1}^\infty$.
\item[{\rm (ii)}]
If $F=(f_n)_{n=1}^\infty$ is a sequence in $\h$   such that $M_{1/m, F, \Phi^\dagger}$is well-defined and
 $$M_{m, \Phi, \Psi}^{-1}  = M_{1/m, F, \Phi^\dagger},$$
then $F$ must be an s-pseudo-dual   of  $\Psi$.
\end{itemize}
\end{prop}

If the symbol is assumed to be bounded, the assumptions that the weighted sequences is a frame is automatically fulfilled, so in summary:

   \begin{cor} \label{ff4} 
Let $\Phi$ and $\Psi$ be frames for $\h$, and let 
the symbol  
$m$ be such that $m_n\neq 0$ for every $n$ and $m\in\ell^\infty$.  
Assume that $M_{m, \Phi, \Psi}$ is invertible. 
Then Propositions \ref{ff3} and  \ref{ff32} apply. 
\end{cor} 
 \bp
Write $M_{m, \Phi,\Psi}$ in the way $M_{(1), m\Phi,\Psi}$. Since $m\Phi$  and $\Psi$ are Bessel sequences and $M_{m, \Phi,\Psi}$ is invertible, \cite[Prop. 3.1]{balsto09new} implies that $m\Phi$ is a frame for $\h$ and thus, Proposition \ref{ff3} applies.
In a similar way, writing $M_{m, \Phi,\Psi}=M_{(1), \Phi,\overline{m}\Psi}$ and using \cite[Prop. 3.1]{balsto09new}, it follows that $\overline{m}\Psi$ is a frame for $\h$ and thus, Proposition \ref{ff32} applies.
\ep

 \v 
  For an illustration of Corollary \ref{ff4}, consider the following example.
	  \begin{ex} {\rm \cite{bstable09}}
Consider 
 
 \ \ \ $\Phi=(e_1,  e_1, - e_1, \   e_2,\ \,    e_2,  \, - e_2,   \, e_3, \  \  e_3, \,  - e_3,  \ldots)$,
 
 \ \ \ $\Psi=(e_1,  e_1,\ \   e_1,\    e_2,   \frac{1}{2}e_2,  \  \frac{1}{2}e_2,  \,  e_3,   \frac{1}{3}e_3,  \ \frac{1}{3}e_3,  \ldots)$,
 
 \ \ \  $m=(\ 1, \ 1, \ \ \ \, 1, \ \ 1, \ \ \, \frac{1}{2}, \ \ \ \ \frac{1}{2}, \ \, 1,\ \ \, \frac{1}{3},\ \ \, \  \frac{1}{3}, \ldots)$.
 
\noindent Then $M_{m,\Phi,\Psi}=Id_{\h}$ and all the assumptions of Corollary \ref{ff4}  are satisfied. 
 \end{ex}

  To illustrate that Proposition \ref{ff3} covers a larger class of multipliers then Corollary \ref{ff4}, consider the following example.
	
  \begin{ex} 
 Consider 
 
 \ \ \ $\Phi=(e_1, e_2, \ \ \frac{1}{2}e_1, e_3,\ \, \frac{1}{2^2}e_1, e_4, \ \ \frac{1}{2^3}e_1, e_5, \ldots)$, 
 
 \ \ \ $\Psi=(e_1, e_2, \frac{1}{\sqrt{2}}e_1, e_3, \frac{1}{\sqrt{2^2}}e_1, e_4, \frac{1}{\sqrt{2^3}}e_1, e_5, \ldots)$,
 
 \ \ \  $m=(\ 1, \ 1,\ \ \ \sqrt{2}, \ 1, \ \ \, \sqrt{2^2}, \ 1, \ \ \, \sqrt{2^3}, \ 1, \ldots)$.

\noindent Then $M_{m,\Phi,\Psi}=Id_{\h}$ and all the assumptions of Proposition \ref{ff3} are satisfied.
 Since $m\notin\ell^\infty$ and $m\Psi$ is not a frame, Corollary \ref{ff4} and  Proposition \ref{ff32} do not apply.

 \end{ex}

\section{ On [Q3]: Inversion using the canonical duals} \label{struct0}

Here we show that in general the answer of $[\bf Q3]$ is negative, providing an example:

\begin{ex}\label{ex3}
Consider the sequences

\vspace{.05in}
\hspace{.367in}
$\Phi=(e_1, e_1, -e_1, e_2, e_2, -e_2, e_3, e_3, -e_3,\ldots),$

\vspace{.05in} \hspace{.365in}
$\Psi=(e_1, e_1, e_1, e_2, e_2, e_2, e_3, e_3, e_3,\ldots),$

\vspace{.05in} \hspace{.365in} $m=(\frac{5+2\sqrt{5}}{5}, \frac{5-2\sqrt{5}}{5},1,
\frac{5+2\sqrt{5}}{5}, \frac{5-2\sqrt{5}}{5},1,
\frac{5+2\sqrt{5}}{5}, \frac{5-2\sqrt{5}}{5},1, \ldots)$, 

\vspace{.05in} \noindent 
for which one has that  $M_{m,\Phi,\Psi}=I_\h$. Here $\widetilde{\Phi}=\frac{1}{3}\Phi$, $\widetilde{\Psi}=\frac{1}{3}\Psi$, and 
$$M_{1/m,\widetilde{\Psi},\widetilde{\Phi}}  = \frac{1}{9}M_{1/m,\Psi,\Phi} 
=I_\h = M^{-1}_{m,\Phi,\Psi}.$$ However, $\Psi$ is not equivalent to $m\Phi$ and 
$\Phi$ is not equivalent to $\overline{m}\Psi$. 

Furthermore,  $M_{m,\Phi,\Psi}$ can not be written as $M_{(1),(c_n\phi_n)_{n=1}^\infty,(d_n\psi_n)_{n=1}^\infty}$ for any $(c_n)_{n=1}^\infty$ and $(d_n)_{n=1}^\infty$ with $(c_n\phi_n)_{n=1}^\infty$ being equivalent to $(d_n\psi_n)_{n=1}^\infty$. Assume the opposite, i.e., that $M_{m,\Phi,\Psi}$ can be written as described. Then all $c_n$ and $d_n$ would be non-zero, 
because $c_n \overline{d_n}=m_n\neq 0$ for all $n$. Denote $k:=\frac{d_1}{c_1}$ 
and let $V:\h\to\h$ denote the bounded bijective operator mapping every $c_n\phi_n$  to the corresponding $d_n\psi_n$. 
Then $c_1 V(e_1) = V(c_1 e_1)= d_1 e_1$ and $-c_3 V(e_1) = V(-c_3 e_1)= d_3 e_1$, implying that  $d_3 = -k c_3$. 
 Then $0<m_1=c_1 \overline{d_1}=\overline{k} \,|c_1|^2$  and $0<m_3 = c_3 \overline{d_3}=- \overline{k} \,|c_3|^2$ leading to a contradiction. 
\end{ex}

The above example actually shows that even the weaker question 

\begin{quote} 
[{$\bf Q3'$}] {\em Under the assumptions of Proposition \ref{q}, is it valid that 
$$
M_{m,\Phi,\Psi}^{-1}=M_{1/m,\widetilde{\Psi},\widetilde{\Phi}} \Leftrightarrow 
 (\Psi \mbox{ is equivalent   to } m\Phi ) \mbox{ or }(\Phi \mbox{ is equivalent   to } \overline{m}\Psi ) ?
$$
 } 
\end{quote}

has a negative answer in general. 

\v
As shown in Proposition \ref{q}, for the class of frame multipliers with constant symbols, the answer of $[Q3]$ is  affirmative. In the following we determine a larger class of multipliers, where the answer of $[Q3]$ is still affirmative.
When a frame $\Phi$ is multiplied with a weight $m$, and the result $m\Phi$ is a frame, i.e. a \emph{weighted frame} \cite{xxljpa1}, the conjugate reciprocal weighted canonical dual frame $\frac{1}{\overline{m}\widetilde{\Phi}} $ is always a dual frame of $m\Phi$. If $\frac{1}{\overline{m}\widetilde{\Phi}} $ is precisely the canonical dual of  $m\Phi$, then the answer of $[Q3]$ is affirmative and this in particular extends the formula (\ref{minvrb}) for Riesz multipliers to a certain class of frame multipliers:

\begin{prop}\label{pr52}
Let a frame $\Phi$ for $\h$ and a sequence $m$ with $m_n\neq 0$ for all $n$ be such that  $m\Phi$ is a frame for $\h$ with the property $\widetilde{m\Phi}=\frac{1}{\overline{m}} \widetilde{\Phi}$. Then for any frame $\Psi$ for $\h$: 
\begin{eqnarray*}
\mbox{$M_{m,\Phi,\Psi}$ is invertible and (\ref{minvrb})  holds} & \Leftrightarrow & \mbox{$\Psi$ is equivalent to $m\Phi$}\\
& \Leftrightarrow & \mbox{$M_{m,\Phi,\Psi}$ is invertible and $\Psi^\dagger = \widetilde{\Psi}$.}
\end{eqnarray*}
\end{prop}
\bp 
 By \cite[Theor. 4.6]{bsreprinv2015}, which concerns frame multipliers with constant symbols, we have that
\begin{eqnarray*}
 \Psi \mbox{ is equivalent to } m\Phi  &\Leftrightarrow&  M_{(1),m\Phi,\Psi} \mbox{ is invertible and } M_{1,m\Phi,\Psi} ^{-1} = M_{(1),\widetilde{\Psi}, \widetilde{m\Phi}} \\
 &\Leftrightarrow& \mbox{$M_{(1),m\Phi,\Psi}$ is invertible and } (M_{(1),m\Phi,\Psi}^{-1} (m_n\phi_n))_{n=1}^\infty=\widetilde{\Psi}. 
\end{eqnarray*} 
Using the assumption  $\widetilde{m\Phi}=\frac{1}{\overline{m}} \widetilde{\Phi}$, it follows that 
$M_{(1),\widetilde{\Psi}, \widetilde{m\Phi}}$ is actually the multiplier $M_{1/m, \widetilde{\Psi}, \widetilde{\Phi}}$ and this completes the proof of the first equivalence. 
For the second equivalence, use  Proposition \ref{ff3}(i).
\ep

\v
In a similar way as in Proposition \ref{pr52}, one can show validity of the following statement:
\begin{prop}\label{pr53}
Let a frame $\Psi$ for $\h$ and a sequence $m$ with $m_n\neq 0$ for all $n$ be such that $\overline{m}\Psi$ is a frame for $\h$ with the property  $\widetilde{\overline{m}\Psi}=\frac{1}{m} \widetilde{\Psi}$. Then for any frame $\Phi$ for $\h$: 
\begin{eqnarray*}
\mbox{$M_{m,\Phi,\Psi}$ is invertible and (\ref{minvrb})  holds} &\Leftrightarrow & \mbox{$\Phi$ is equivalent to $\overline{m}\Psi$}\\
&\Leftrightarrow &\mbox{$M_{m,\Phi,\Psi}$ is invertible and $\Phi^\dagger = \widetilde{\Phi}$.}
\end{eqnarray*}
\end{prop}

Finally we can state a result for symbols with constant amplitude but varying phase:

\begin{cor}\label{cor54}
Let $\Phi$ and $\Psi$ be frames for $\h$ and let $m$ be such that $|m_n|=c\neq 0$ for every $n$. Then 
\begin{eqnarray*}
M_{m,\Phi,\Psi} \mbox{ is invertible and (\ref{minvrb})  holds } & \Leftrightarrow & \Psi \mbox{ is equivalent to } m\Phi \\
 & \Leftrightarrow & \Phi \mbox{ is equivalent to } \overline{m}\Psi.
 \end{eqnarray*}
\end{cor}
\bp  Clearly, $m\Phi$ and $\overline{m}\Psi$ are frames for $\h$. Furthermore, 
 $S_{m\Phi} = c^2 S_\Phi $,  
implying that 
$$S_{m\Phi}^{-1}(m_n \phi_n) =
 \frac{m_n}{c^2}\,  S_\Phi^{-1} (\phi_n)= \frac{1}{\overline{m_n}} \, S_\Phi^{-1} (\phi_n).$$
Therefore,  $\widetilde{m\Phi}=\frac{1}{\overline{m}} \widetilde{\Phi}$. In a similar way, it follows that $\widetilde{\overline{m}\Psi}=\frac{1}{m} \widetilde{\Psi}$. Now Propositions \ref{pr52} and \ref{pr53} complete the proof. 
\ep 

\sloppy
\v As a simple example illustrating Corollary \ref{cor54}, consider the frames 
$\Phi=(e_1, e_1, e_2, e_2, e_3, e_3, \ldots)$ and $\Psi=(e_1, -e_1, e_2, -e_2, e_3, -e_3, \ldots)$, and the sequence $m=(1,-1,1,-1,1,-1,\ldots)$, for which one has $\Psi=m\Phi$ and $\Phi=\overline{m}\Psi$.

\vspace{13pt}
\noindent {ACKNOWLEDGEMENT:} The authors acknowledge support from the Austrian Science Fund (FWF) START-project FLAME ('Frames and Linear Operators for Acoustical Modeling and Parameter Estimation'; Y 551-N13).

\end{document}